\documentclass[conference]{IEEEtran}
\usepackage[utf8]{inputenc}
\usepackage{amssymb, amsmath, graphicx, subfigure, enumerate, listings,multicol,amsfonts}
\usepackage{amsthm,alltt}
\usepackage{graphicx,ctable,booktabs}
\usepackage{mathtools}
\usepackage[boxed]{algorithm2e}
\usepackage{mathdots}
\usepackage{fancyhdr} 
\usepackage{cleveref}
\usepackage{amssymb}
\usepackage{amsthm}
\usepackage{amscd}
\usepackage{amscd}
\usepackage{enumerate}

\makeatletter

\newtheorem{theorem}{Theorem}[section]

\newtheorem{lemma}[theorem]{Lemma}
\newtheorem{conjecture}[theorem]{Conjecture}

\newtheorem{proposition}[theorem]{Proposition}

\newcommand{\C}{\mathbb{C}}
\newcommand{\R}{\mathbb{R}}

\begin{document}

\title{Frame Potentials and Orthogonal Vectors}

\author{
\IEEEauthorblockN{Josiah {\sc Park}}

\IEEEauthorblockA{School of Mathematics \\ Georgia Tech \\ Atlanta, GA, 30332 \\\
$${\tt j.park@gatech.edu}}}

\maketitle

\begin{abstract}
An extension is given of a recent result of Glazyrin, showing that an orthonormal basis $\{e_{i}\}_{i=1}^{d}$ joined with the vectors $\{e_{j}\}_{j=1}^{m}$, where $1\leq m < d$ minimizes the $p$-frame potential for $p\in[1,2\log{\frac{2m+1}{2m}}/\log{\frac{m+1}{m}}]$ over all collections of $N=d+m$ vectors $\{x_1,\dots,x_N \}$ in $\mathbb{S}^{d-1}$.
\end{abstract}

\section{Introduction} \label{intro}
For a set of unit vectors $X=\{x_1,\dots,x_N\}\subset \mathbb{S}^{d-1}$, an interesting quantity associated with $X$ is the $p$-frame energy $\sum\limits_{i\neq j} |\langle x_{i},x_{j} \rangle|^p.$ This energy perhaps appeared earliest for even $p$ and vectors in $\mathbb{S}^4$ in Hilbert's 1909 solution to Waring's problem \cite{DH}. For even $p$ the energy has close ties to objects called spherical $t$-designs, certain configurations which act as nodes for integration over the sphere, appearing for instance in a 1981 paper of Goethals and Seidel \cite{GS}. More recently, the terms $p$-frame energy or potential are used sometimes interchangeably to describe the $\ell_p$ norm of the off-diagonal elements in the Gram matrix of a collection of unit vectors, this term originating in a paper of Benedetto and Fickus \cite{BF}. In their 2003 paper they introduced this term for the $p=2$ energy after observing that minimizers are precisely what are known as finite unit norm tight frames (FUNTFs).

\indent For even $p$, minimizers of this quantity for sufficiently many points on the sphere are $t$-designs. Associated identities holding more generally for weighted designs show also that some minimizers have interpretations as minimal isometric embeddings of finite dimensional $\ell_2$ spaces into higher dimensional $\ell_p$ spaces \cite{OS}. In projective space, the analogous minimizers for $p$ even (and unit norm vectors in $\C^{d}$) are known to be projective $t$-designs and have optimal properties for measuring quantum states \cite{AR}.

\indent Describing minimizers for the $p$-frame potential for $p$ not even appears to be a difficult problem, and in general not much is known about the structure of minimizers outside a few exceptional cases (some results can found in the papers of Ehler and Okoudjou in this line \cite{EO} or in the recent pre-print \cite{CGEKO}). A large part of the literature surrounding these energies focuses on their relationship  with certain symmetric minimal coherence systems of vectors known as equiangular tight frames (ETFs), studies of which first appeared in the discrete geometry community \cite{LS}. An elementary argument shows that ETFs minimize the $p$-frame energy for $p>2$, when these systems exist \cite{OO}. It a well-known open problem to determine when an ETF of $N$ unit vectors exists generally and much of the evidence for existence of ETFs outside of the real case is due to the observation that they are minimizers of a range of energies \cite{CKM}.

In this note it is demonstrated that the method developed by Glazyrin in \cite{AG} for describing minimizers of $p$-frame energies has further applications. Adopting the notation used there for the Gram matrix $A=X'X$ of a system of unit vectors $X$, the $p$-frame energy may be given alternatively by $E_{p}(A)=\sum\limits_{i\neq j} |A_{i,j}|^p $. The main observation here is the following result.

\begin{theorem} For $p\in[1,2\log{\frac{2m+1}{2m}}/\log{\frac{m+1}{m}}]$, $1\leq m< d$ and real $(d+m)\times (d+m)$ matrix $A$ of rank $d$ with ones along the diagonal, 
$$E_{p}(A)=\sum\limits_{i\neq j}|A_{i,j}|^{p} \geq 2m.$$
\end{theorem}
 The proof for the above inequality is an extension of the method used in \cite{AG}, and by restriction to $m=1$, one obtains the result proved there. In the limit, it is known that the unique symmetric Borel probability measure which minimizes the $p$-frame energy on the sphere $\mathbb{S}^{d-1}$ equally distributes mass over the vertices of a cross-polytope whenever $p\in(0,2)$ \cite{EO}. These energies do not depend on the sign of any vector and so one can reflect any vector about the origin to obtain the same energy. For this reason it makes the most sense to consider the energy projectively, that is with vectors constrained to lie in one hemisphere. With this in mind, the above shows that for $N=d+m$ vectors, $1\leq m<d$, and $p$ in a certain range near $p=1$, the support for the finite minimization problem agrees with that of the limiting distribution.

 It will be necessary to introduce a related optimization problem to minimizing $E_{p}$ found in the previously mentioned reference \cite{AG} in order to state the relevant steps in the proof of optimality of the orthonormal sequence for the above mentioned range of $p$.

\section{Repeated Ortho-sequence Minimizes $E_{p}(A)$} \label{proofsection}

 Define $f_{c,p}(t)=\left(\frac{t}{c-t}\right)^{\frac{p}{2}}$ and set $M(c,p,N)$ to be the optimal value in the optimization problem
\[ \min\left\{\sum\limits_{i=1}^{N} f_{c,p}(t_{i})\ |\ \sum\limits_{i=1}^{N} t_{i}=1,\ t_{i}\in [0,c) \right\}. \]
The following inequality for $E_p$ and $M(c,p,N)$ is proved in \cite[Lemma 2.2]{AG}.

\begin{proposition} For any real $N\times N$ matrix $A$ of rank $d$ with unit diagonal elements, 
\[E_{p}(A) \geq M\left(\frac{1}{N-d},p,N \right),\ \text{for}\  1\leq p \leq 2. \]

\end{proposition}

\indent By the above proposition, in order to prove the theorem it suffices to show $M(\frac{1}{m},p,N)\geq 2m$. The following observation, used in the proof of the case $m=1$ in \cite{AG}, will be applied below (which is obtained by use of concavity/convexity of $f$ and Jensen's and Karamata's inequality):

\begin{lemma} Set $\alpha=\frac{1}{2}-\frac{p}{4}$. For $p\in[1,2]$, $M(\frac{1}{m},p,N)$ is minimized for $t_{j}$ of the form
\begin{enumerate}
\item[(i)] $t_1=\dots=t_{k}=\frac{1}{k},\ t_{k+1}=\dots=t_{n}=0$, where $\frac{1}{k}\geq \alpha$ 
\item[{}] \hspace{3.5 cm} or 
\item[(ii)] $t_1=\dots=t_k=x$, $t_{k+1}=1-kx$,\ \\ $t_{k+2}=\dots=t_{N}=0$, where $x\geq \alpha$, $0<1-kx<\alpha$.
\end{enumerate}
\end{lemma}
The proof of the main theorem is now given.
\begin{proof} 
\indent Set $p_{0,m}=2\log{\frac{2m+1}{2m}}/\log{\frac{m+1}{m}}$ and $q_{m}=\frac{p_{0,m}}{2}$. Consider the first case in the above lemma, $t_1=\dots=t_{k}=\frac{1}{k},\ t_{k+1}=\dots=t_{n}=0$, where $\frac{1}{k}\geq \alpha$. In this case, for $p<p_{0,m}$ $kf_{1,p}(\frac{1}{k})=\frac{k}{(k-1)^{\frac{p}{2}}}$ takes minimal value $2m$.

\indent In the second case, $x<\frac{1}{k}$ and $x\geq \alpha\geq \frac{1}{2}-\frac{p_{0,m}}{4}$ so that $k$ can take (integer) values only in $[m,4m]$. To show $E_{p}\geq 2m$ for $p\leq p_{0,m}$, it suffices then to show for all $m\leq j\leq 4m$, and all $x$ in $I=(\frac{1}{j+1},\frac{1}{j})$ that 
\begin{align*}    g_{j}(x) =  j \left( \frac{mx}{1-mx} \right)^{q_{m}}+\left(\frac{m(1-jx)}{1-m(1-jx)}\right)^{q_{m}}  \end{align*} 
satisfies $g_{j}(x)\geq 2m$. This will be demonstrated using properties specific to $g_{j}(x)$, namely that the function has at most one critical point, $g_{j}^{\prime}(x)=0$, inside the interval $I$. Taking derivatives,
 \begin{align*} g^{\prime}_{j}(x)&=q_{m} jm  \left( \frac{mx}{1-mx}  \right)^{q_{m}-1} \frac{1}{(1-mx)^2} \\
 & -q_{m} jm \left(\frac{m(1-jx)}{1-m(1-jx)}\right)^{q_{m}-1} \frac{1}{(1+m(-1+jx))^2}, \end{align*}

 so that $g'_{j}(x)=0$  gives
 \begin{align*}  \left( \frac{x(1+m(-1+jx))}{(1-mx)(1-jx)} \right)^{q_{m}-1} &= \frac{(1-mx)^2}{(1+m(-1+jx))^2} \\ 
\left( \frac{x(1+m(-1+jx))}{(1-mx)(1-jx)} \right)^{q_{m}+1} &=  \frac{x^2}{(1-jx)^2}  \\
  \frac{1+m(-1+jx)}{1-mx}  &= \left( \frac{x}{(1-jx)} \right)^{\frac{2}{q_{m}+1}-1} \\
 \frac{1-mx}{1+m(-1+jx)}  &= \left( \frac{x}{(1-jx)} \right)^{1-\frac{2}{q_{m}+1}}.
\end{align*} %
\vspace{1 mm}
Calling the function on the left in the above expression $f(x)$ and the function on the right $g(x)$,
\begin{align*}
f^{\prime \prime}(x) = \frac{2j(1+j-m)m^2}{(1+m(-1+jx))^3}>0 \text{ on }I, \end{align*} while letting $\alpha=1-\frac{2}{q_{m}+1}$, \begin{align*} g^{\prime \prime}(x)  = \frac{\alpha (\frac{x}{1-jx})^\alpha (-1+\alpha+2jx)}{x^2(jx-1)^2}<0 \text{ on } I, \end{align*}   since $\alpha<0$. Thus $f(x)$ is convex on $I$, while $g(x)$ is concave on $I$. Since $f(\frac{1}{j+1})=g(\frac{1}{j+1})$ and $f'(\frac{1}{j+1})\leq g'(\frac{1}{j+1})$ when $j<4m$ it must be the case then that $f(x)=g(x)$ for exactly one point $x\in I$, ($x\neq \frac{1}{j+1},\frac{1}{j}$). Note that when $j=4m$ there are no critical points in $I$. Now, \begin{align*} g_{j}'\left(\frac{1}{j+1}\right)=0\ \text{ and } \lim\limits_{x\rightarrow \frac{1}{j}} g_{j}'\left(x\right)=-\infty. \end{align*} \indent So the critical points then correspond to local maxima of $g_{j}(x)$ and it suffices to check the value of $g_{j}(x)$ at the endpoints in $I$ for each $m\leq j\leq 4m$ to establish the desired lower bound. These values are
\begin{align*} g_{j}\left(\frac{1}{j+1}\right)&=(1+j)\left(\frac{m}{1+j-m}\right)^{q_m}, \\ \  g_{j}\left(\frac{1}{j}\right)&=j\left(\frac{m}{j-m}\right)^{q_m}. \end{align*}
\indent Each value may be checked to be greater or equal to $2m$ by minimizing with respect to $j$. Taking a derivative in $j$ gives a decreasing then increasing expression with a zero between $2m$ and $2m+1$. These closest values then minimize the expression over all feasible positive integers $j\geq m$ and the minimal value for all cases is $2m$. \end{proof}
\vspace{-8 mm} \section{Discussion} \label{furtherdiscussion} 
As was noted in \cite{AG}, the above argument applies to the problem of minimizing $E_{p}(A)$ over $N\times N$ matrices over $\mathbb{F}=\mathbb{R},\mathbb{C}$, or $\mathbb{H}$, real and complex numbers or quaternions. For $N=d+1$ the range of $p$ for which the orthogonal construction above is expected to be optimal for $E_{p}$ is $p\in[0,\frac{\log{3}}{\log{2}}]$ and this question is part of a more general conjecture by Chen, Gonzales, Goodman, Kang, and Okoudjou \cite{CGEKO} about minimizers of $E_{p}(A)$ with $p\in[0,2]$ (and $N=d+1$). As was noted also in \cite{AG}, the bound given by the main theorem here does not extend fully to this conjectured range. How far from sharp the above bound is for $m>1$ appears to be an interesting question. 

We briefly look into this question now, building on some recent observations from \cite{CGEKO}. It was suggested from a numerical study that for $N=5$ points on the unit circle there may be a transition around $p=1.78$ for which the frame energy changes from being minimized on $\{e_1,e_1,e_2,e_2,e_2\}$ to a configuration of the form $\{x,x,y,y,z\}$. One example of a Gram matrix from a system of vectors which can take this form, (but need not generally) is the matrix $$A=\left( \begin{matrix} \ 1 & \ 1 & \ 0 & \ \alpha & -\alpha \\ \ 1 & \ 1 & \ 0 & \ \alpha & -\alpha \\ \ 0 & \ 0 & \ 1 & \sqrt{1-\alpha^2} & \sqrt{1-\alpha^2} \\ \ \alpha & \ \alpha & \sqrt{1-\alpha^2} & \ 1 & \ \beta \\ -\alpha & -\alpha & \sqrt{1-\alpha^2} & \ \beta & \ 1  \end{matrix} \right).$$
 Since $A$ is a rank-two matrix, $$\det \left( \begin{matrix} \ 1 & \ \alpha & -\alpha \\ \ \alpha & \ 1 & \ \beta \\ -\alpha & \ \beta & \ 1 \end{matrix} \right) =0,$$ so that $\beta=-1$ or $\beta=1-2\alpha^2$. The first value gives a larger $E_p$ value, so suppose instead that $\beta=1-2\alpha^2$. Then for this $A$,
 $${\small E_p(A)=2+8\alpha^{p}+2(1-2\alpha^2)^p+(1-\alpha^2)^\frac{p}{2}} \hspace{1 mm} \text{ and }$$
\small \begin{eqnarray*} \frac{dE_p(A)}{d \alpha}=p(8\alpha^{p-1}-8\alpha(1-2\alpha^2)^{p-1}-4\alpha(1-\alpha^2)^{\frac{p}{2}-1}).\end{eqnarray*} \normalsize
Note now that the value of $E_p$ on the repeated orthonormal sequence $\{e_1,e_1,e_2,e_2,e_2\}$ is $4$.
It remains now to consider solutions $(\alpha,p)$ to the system 
 \begin{eqnarray*}  8\alpha^{p}+2(1-2\alpha^2)^p+(1-\alpha^2)^{\frac{p}{2}}=4, \\ \hspace{2 mm} p(8\alpha^{p-1}-8\alpha(1-2\alpha^2)^{p-1}-4\alpha(1-\alpha^2)^{\frac{p}{2}-1})=0.  \end{eqnarray*} 
Given that one may not expect such a system to have solutions necessarily expressible via elementary functions, looking numerically for a solution gives the values of $\alpha$ and $p$ below
\begin{align*} \alpha&=&0.43421690071432109168188584186122094, \\
 p&=&1.77766251887018589539510545748522601. \end{align*}
\indent Replacing $4$ on the right hand side of the first equation above with $4$ minus a small quantity and repeating the root finding procedure provides a pairing $(\alpha,p)$ with a smaller corresponding value of $E_p$ than $4$ (which can be checked to be valid by truncating the numerical solution at a given precision, noting that this $\alpha$ will still be feasible). 

\indent After experimenting numerically it appears one can extend these observations similarly to the case of $N=7$, where the transition value there appears to be about $p=1.840321171266$. In both of the above cases, observations only provide evidence that the threshold can occur no later than the $p$ value above. A more general picture is suggested by further experiments.

\begin{conjecture} Let $N=m+kd$ points be given in $\mathbb{S}^{d-1}$, with $1\leq m<d$, $d\geq 2$, and gram matrix $A\in\R^{N\times N}$. Then there is a value of $p_0$, independent of dimension $d$ and excess $m$, such that the repeated orthonormal sequence $\{e_{j \mod d} \}_{j=1}^N$ minimizes $E_p$ over all size $N$ systems of unit vectors (with value $E_p(A)=d(k^2-k)+2k$) for $p<p_0$ and the minimum value of $E_p(A)$  satisfies $E_p(A)<d(k^2-k)+2k$ when $p>p_0$. Further $p_0=p_0(k)$ satisfies $p_0(k)\rightarrow 2$ as $k\rightarrow \infty$.
\end{conjecture}

\section{Acknowledgements} \label{acknowledge}
 The author was supported in part by grant from the US National Science Foundation, DMS-1600693. The author would like to acknowledge fruitful discussions with Alexey Glazyrin. The author would also like to acknowledge helpful discussions about related topics with Dmitriy Bilyk, Ryan Matzke, and Oleksandr Vlasiuk in connection with a forthcoming paper.

\end{document}